# Partial-State DADS Control for Matched Unmodeled Dynamics

**Iasson Karafyllis[*] and Miroslav Krstic[**]**

[*]Dept. of Mathematics, National Technical University of Athens, Zografou Campus, 15780, Athens, Greece,
email: iasonkar@central.ntua.gr; iasonkaraf@gmail.com

[**]Dept. of Mechanical and Aerospace Eng., University of California, San Diego, La Jolla, CA 92093-0411, U.S.A., email: krstic@ucsd.edu

## Abstract

We extend the Deadzone-Adapted Disturbance Suppression (DADS) control to time-invariant systems with dynamic uncertainties that satisfy the matching condition and for which no bounds for the disturbance and the unknown parameters are known. This problem is equivalent to partial-state adaptive feedback, where the states modeling the dynamic uncertainty are unmeasured. We show that the DADS controller can bypass small-gain conditions and achieve robust regulation for systems in spite of the fact that the strength of the interconnections has no known bound. Moreover, no gain and state drift arise, regardless of the size of the disturbances and unknown parameters. Finally, the paper provides the detailed analysis of a control system where the unmeasured state (or the dynamic uncertainty) is infinite-dimensional and described by a reaction-diffusion Partial Differential Equation, where the diffusion coefficient and the reaction term are unknown. It is shown that even in the infinite-dimensional case, a DADS controller can be designed and guarantees robust regulation of the plant state.

## 1. Introduction

The issue of robustness in adaptive control is very important and has been studied by many researchers in control/stability theory. When dealing with time-invariant nonlinear control systems and no persistence of excitation condition is assumed, numerous approaches have been proposed in the literature: leakage (see [2, 3, 20, 23]), nonlinear damping (see [14, 15, 7, 8]), dynamic (high) gains or gain adjustment (see [4, 13, 17, 18, 22]), projection methodologies (see [2] and Appendix E in [14]), supervision for direct adaptive schemes (see [1]) and deadzone in the update law (introduced in the paper [21] and well explained in the book [2]).

Recently, a new approach was proposed in the literature (see [9, 10, 11]): the Deadzone-Adapted Disturbance Suppression (DADS) controller. DADS provides a simple, direct, adaptive control scheme that combines three elements: (a) nonlinear damping (as in [14, 8]), (b) single-gain adjustment (the dynamic feedback has only one state), and (c) deadzone in the update law. The proposed adaptive control scheme achieves attenuation of the plant state to an assignable small level, despite the presence of (time-varying) disturbances and unknown (time-varying) parameters of arbitrary and unknown bounds. Moreover, the DADS controller prevents gain and state drift regardless of the size of the disturbance and unknown parameter. The latter property is a consequence of the practical Uniform Bounded-Input-Bounded-State (p-UBIBS) property and the practical Input-to-Output Stability (p-IOS) property that both hold for the closed-loop system.

In the present note, we extend the DADS controller to time-invariant systems with dynamic uncertainties that satisfy the so-called matching condition (see [14] for the definition of the matching condition) for which no bounds for the disturbance and the unknown parameters are known. This specific control problem can equivalently be posed as a partial-state adaptive feedback design problem, where the states describing the effect of the dynamic uncertainty are the unmeasured states. We show that the DADS control can allow the "escape" from small-gain conditions (proposed in [5]) and achieve robust regulation for systems where the "strength" of the interconnections has no known bound. Moreover, as in [9, 10, 11], no gain and state drift phenomena are present, regardless of the size of the disturbances and unknown parameters.

The present work goes one step beyond the standard analysis of dynamic uncertainties for adaptive controllers. Section 3 provides the detailed analysis of a control system where the unmeasured state (or the dynamic uncertainty) is infinite-dimensional. It is shown that the DADS controller can be applied with no problem and guarantee robust regulation of the plant state despite the presence of (time-varying) disturbances and unknown (time-varying) parameters of arbitrary and unknown bounds. Therefore, DADS is not limited to the finite-dimensional case and can deal with certain infinite-dimensional cases as well.

The structure of the paper is as follows. First, we start with a subsection that provides the notation and all the stability notions used in the paper (some of them are modifications of well-known notions presented in [6, 7, 11, 12, 24]). All the main results are stated and discussed in Section 2. Section 3 is devoted to the study of the case where the dynamic uncertainty is infinite-dimensional and is described by a reaction-diffusion Partial Differential Equation (PDE), where the diffusion coefficient and the reaction term are unknown. Section 4 of the paper contains the proofs of all main results. The concluding remarks of the present work are provided in Section 5.

**Notation and Basic Notions.** Throughout this paper, we adopt the following notation.

* $\mathbb{R}_+ := [0, +\infty)$. For a vector $x \in \mathbb{R}^n$, $|x|$ denotes its Euclidean norm and $x'$ denotes its transpose. We use the notation $x^+$ for the positive part of the real number $x \in \mathbb{R}$, i.e., $x^+ = \max(x, 0)$.

* Let $D \subseteq \mathbb{R}^n$ be an open set and let $S \subseteq \mathbb{R}^n$ be a set that satisfies $D \subseteq S \subseteq cl(D)$, where $cl(D)$ is the closure of $D$. By $C^0(S; \Omega)$, we denote the class of continuous functions on $S$, which take values in $\Omega \subseteq \mathbb{R}^m$. By $C^k(S; \Omega)$, where $k \geq 1$ is an integer, we denote the class of functions on $S \subseteq \mathbb{R}^n$, which take values in $\Omega \subseteq \mathbb{R}^m$ and have continuous derivatives of order $k$. In other words, the functions of class $C^k(S; \Omega)$ are the functions which have continuous derivatives of order $k$ in $D = \text{int}(S)$ that can be continued continuously to all points in $\partial D \cap S$. When $\Omega = \mathbb{R}$ then we write $C^0(S)$ or $C^k(S)$. A function $f \in \bigcap_{k=0}^{\infty} C^k(S; \Omega)$ is called a smooth function.

* By $L^\infty(\mathbb{R}_+;\Omega)$, where $\Omega \subseteq \mathbb{R}^n$, we denote the class of essentially bounded, Lebesgue measurable functions $d:\mathbb{R}_+ \to \Omega$. For $d \in L^\infty(\mathbb{R}_+;\Omega)$ we define $\|d\|_\infty = \sup_{t \geq 0}(|d(t)|)$, where $\sup_{t \geq 0}(|d(t)|)$ is the essential supremum. When $\Omega = \mathbb{R}$ then we simply write $L^\infty(\mathbb{R}_+)$ instead of $L^\infty(\mathbb{R}_+;\mathbb{R})$.

* By $K$ we denote the class of increasing continuous functions $a:\mathbb{R}_+ \to \mathbb{R}_+$ with $a(0)=0$. By $K_\infty$ we denote the class of increasing continuous functions $a:\mathbb{R}_+ \to \mathbb{R}_+$ with $a(0)=0$ and $\lim_{s \to +\infty}(a(s)) = +\infty$. By $KL$ we denote the set of all continuous functions $\sigma:\mathbb{R}_+ \times \mathbb{R}_+ \to \mathbb{R}_+$ with the properties: (i) for each $t \geq 0$ the mapping $\sigma(\cdot,t)$ is of class $K$; (ii) for each $s \geq 0$, the mapping $\sigma(s,\cdot)$ is non-increasing with $\lim_{t \to +\infty}(\sigma(s,t)) = 0$.

* Let $S \subseteq \mathbb{R}^n$ be a non-empty set with $0 \in S$. We say that a function $V:S \to \mathbb{R}_+$ is positive definite if $V(x) > 0$ for all $x \in S$ with $x \neq 0$ and $V(0)=0$. We say that a continuous function $V:S \to \mathbb{R}_+$ is radially unbounded if the following property holds: "for every $M > 0$ the set $\{x \in S : V(x) \leq M\}$ is compact". For $V \in C^1(S;\mathbb{R}_+)$ we define $\nabla V(x) = \left(\frac{\partial V}{\partial x_1}(x),...,\frac{\partial V}{\partial x_n}(x)\right)$.

* Let $a < b$ be given constants. $L^2(a,b)$ is the set of equivalence classes of Lebesgue measurable functions $u:(a,b) \to \mathfrak{R}$ with $\|u\| := \left(\int_a^b |u(x)|^2 dx\right)^{1/2} < +\infty$. For an integer $k \geq 1$, $H^k(a,b)$ denotes the Sobolev space of functions in $L^2(a,b)$ with all its weak derivatives up to order $k \geq 1$ in $L^2(a,b)$. The closure of the class of functions $f \in C^1((a,b))$ with compact support in $H^1(a,b)$ is denoted by $H_0^1(a,b)$.

We next recall certain notions of output stability that are used in this work (see also [11] and references therein). Let $D \subseteq \mathbb{R}^p$ be a given closed set with $0 \in D$, $f:\mathbb{R}^n \times D \to \mathbb{R}^n$ be a locally Lipschitz with respect to $x \in \mathbb{R}^n$ mapping with $f(0,0)=0$ and $h:\mathbb{R}^n \to \mathbb{R}^p$ be a continuous mapping with $h(0)=0$. Consider the control system

$$\dot{x} = f(x,d), \ x \in \mathbb{R}^n, \ d \in D \tag{1.1}$$

with output

$$y = h(x) \tag{1.2}$$

We assume that system (1.1) is forward complete, i.e., for every $x_0 \in \mathbb{R}^n$ and for every Lebesgue measurable and locally essentially bounded input $d:\mathbb{R}_+ \to D$ the unique solution $x(t) = \phi(t,x_0;d)$ of the initial-value problem (1.1) with initial condition $x(0) = x_0$ corresponding to input $d:\mathbb{R}_+ \to D$ exists for all $t \geq 0$. We use the notation $y(t,x_0;d) = h(\phi(t,x_0;d))$ for all $t \geq 0$, $x_0 \in \mathbb{R}^n$ and for every Lebesgue measurable and locally essentially bounded input $d:\mathbb{R}_+ \to D$.

We say that system (1.1), (1.2) is *Input-to-Output Stable (IOS)* if there exist a function $\beta \in KL$ and a non-decreasing, continuous function $\gamma : \mathbb{R}_+ \to \mathbb{R}_+$ with $\gamma(0) = 0$ such that the following estimate holds for all $x_0 \in \mathbb{R}^n$, $t \geq 0$ and for every $d \in L^\infty(\mathbb{R}_+; D)$:

$$|y(t, x_0; d)| \leq \beta(|x_0|, t) + \gamma(\|d\|_\infty) \tag{1.3}$$

We say that system (1.1), (1.2) is *practically Input-to-Output Stable (p-IOS)* if there exist a function $\beta \in KL$, a non-decreasing, continuous function $\gamma : \mathbb{R}_+ \to \mathbb{R}_+$ with $\gamma(0) = 0$ and a constant $\alpha > 0$ such that the following estimate holds for all $x_0 \in \mathbb{R}^n$, $t \geq 0$ and for every $d \in L^\infty(\mathbb{R}_+; D)$:

$$|y(t, x_0; d)| \leq \beta(|x_0|, t) + \gamma(\|d\|_\infty) + \alpha \tag{1.4}$$

The constant $\alpha > 0$ is called the *residual constant* while the function $\gamma$ is called the *gain function of the input $d \in D$ to the output $y$*.

We say that system (1.1), (1.2) satisfies the *practical Output Asymptotic Gain (p-OAG)* property if there exists a non-decreasing, continuous function $\tilde{\gamma} : \mathbb{R}_+ \to \mathbb{R}_+$ with $\tilde{\gamma}(0) = 0$ and a constant $\tilde{\alpha} > 0$ such that the following estimate holds for all $x_0 \in \mathbb{R}^n$ and for every $d \in L^\infty(\mathbb{R}_+; D)$:

$$\limsup_{t \to +\infty} (|y(t, x_0; d)|) \leq \tilde{\gamma}(\|d\|_\infty) + \tilde{\alpha} \tag{1.5}$$

The constant $\tilde{\alpha} > 0$ is called the *asymptotic residual constant* while the non-decreasing, continuous function $\tilde{\gamma}$ is called the *asymptotic gain function of the input $d \in D$ to the output $y$*.

When $\tilde{\gamma} \equiv 0$ we say that system (1.1), (1.2) satisfies the *zero practical Output Asymptotic Gain property (zero p-OAG)*.

When $h(x) = x$ then the word "output" in the above properties is either replaced by the word "state" (e.g., ISS, p-ISS) or is omitted (e.g., p-AG, zero p-AG).

We say that system (1.1) satisfies the *practical Uniform Bounded-Input-Bounded-State (p-UBIBS)* property if there exists a function $\bar{\gamma} \in K_\infty$ and a constant $\bar{\alpha} > 0$ such that the following estimate holds for all $x_0 \in \mathbb{R}^n$ and for every $d \in L^\infty(\mathbb{R}_+; D)$:

$$\sup_{t \geq 0}(|\phi(t, x_0; d)|) \leq \bar{\gamma}(|x_0|) + \bar{\gamma}(\|d\|_\infty) + \bar{\alpha} \tag{1.6}$$

The p-UBIBS property with $\bar{\alpha} = 0$ is called the UBIBS property. Clearly, the p-UBIBS property is equivalent to the existence of a continuous function $B : \mathbb{R}^n \times \mathbb{R}_+ \to \mathbb{R}_+$ for which the following estimate holds for all $x_0 \in \mathbb{R}^n$ and for every $d \in L^\infty(\mathbb{R}_+; D)$:

$$\sup_{t \geq 0}(|\phi(t, x_0; d)|) \leq B(x_0, \|d\|_\infty) \tag{1.7}$$

## 2. Main Results

In this work we study nonlinear control systems of the form

$$\dot{y} = f(y) + g(y)\left(u + \varphi'(y,w)\theta + A'(y,w)d\right) \quad (2.1)$$
$$y \in \mathbb{R}^n, u \in \mathbb{R}, d \in \mathbb{R}^m, \theta \in \mathbb{R}^p$$

$$\dot{w} = h(y,w,\delta) \quad (2.2)$$
$$w \in \mathbb{R}^l, \delta \in \mathbb{R}^q$$

where $f, g : \mathbb{R}^n \to \mathbb{R}^n$, $\varphi : \mathbb{R}^n \times \mathbb{R}^l \to \mathbb{R}^p$, $A : \mathbb{R}^n \times \mathbb{R}^l \to \mathbb{R}^m$, $h : \mathbb{R}^n \times \mathbb{R}^l \times \mathbb{R}^q \to \mathbb{R}^l$ are smooth mappings with $f(0) = 0$, $\varphi(0,0) = 0$, $h(0,0,0) = 0$, $(y,w) \in \mathbb{R}^n \times \mathbb{R}^l$ is the plant state, $u \in \mathbb{R}$ is the control input and $\theta \in \mathbb{R}^p$, $d \in \mathbb{R}^m$, $\delta \in \mathbb{R}^q$ are unknown disturbances. The reader should notice the difference between $\theta$ and $(d,\delta)$ when both are considered to be perturbations: while $\theta \in \mathbb{R}^p$ is a vanishing perturbation (due to the fact that $\varphi(0,0) = 0$), $(d,\delta) \in \mathbb{R}^m \times \mathbb{R}^q$ is -in general- a non-vanishing perturbation. Systems of the form (2.1), (2.2) are systems that satisfy a sort of a matching condition, i.e., the effect of both $\theta \in \mathbb{R}^p$ and $d \in \mathbb{R}^m$ can be cancelled by the control input $u \in \mathbb{R}$ if they are known. We assume next that $d \in L^\infty(\mathbb{R}_+; \mathbb{R}^m)$, $\theta \in L^\infty(\mathbb{R}_+; \mathbb{R}^p)$, $\delta \in L^\infty(\mathbb{R}_+; \mathbb{R}^q)$ but we assume no bounds for $\theta \in \mathbb{R}^p$, $\delta \in \mathbb{R}^q$ and $d \in \mathbb{R}^m$.

Our aim is to design a partial-state robust adaptive feedback law that depends only on $y \in \mathbb{R}^n$.

To this purpose we employ the following assumption.

**Assumption (A):** *There exist smooth mappings* $V, Q : \mathbb{R}^n \to \mathbb{R}_+$, $k : \mathbb{R}^n \to \mathbb{R}$, $\Phi, R : \mathbb{R}^l \to \mathbb{R}_+$, $\mu : \mathbb{R}^n \to (0, +\infty)$ *with* $k(0) = 0$, $V, Q, R, \Phi$ *being positive definite and radially unbounded, constants* $r, \Lambda > 0$ *and a non-decreasing function* $\gamma \in C^0(\mathbb{R}_+; \mathbb{R}_+)$ *such that the following inequalities hold for all* $y \in \mathbb{R}^n$, $w \in \mathbb{R}^l$ *and* $\delta \in \mathbb{R}^q$:

$$\nabla V(y)\left(f(y) + g(y)k(y)\right) \leq -rQ(y) \quad (2.3)$$

$$\nabla \Phi(w) h(y,w,\delta) \leq -R(w) + \gamma(|\delta|) + Q(y) \quad (2.4)$$

$$|\varphi(y,w) - \varphi(y,0)|^2 + |A(y,w) - A(y,0)|^2 \leq \mu(y)(R(w) + \Lambda) \quad (2.5)$$

$$|\varphi(y,0)|^2 \leq \mu(y)\left(Q(y) + \Lambda + (\nabla V(y)g(y))^2\right) \quad (2.6)$$

Assumption (A) guarantees:

(i) that the feedback law $u = k(y)$ is a global feedback stabilizer for the subsystem (2.1) when $\theta = 0$ and $d = 0$; this is a consequence of (2.3),

(ii) that the feedback law $u = k(y)$ achieves the p-ISS property for system (2.1), (2.2) when $\theta = 0$ and $d = 0$ with $\delta \in \mathbb{R}^q$ as input; this is a consequence of (2.3), (2.4) and the use of a Lyapunov function of the form $U(y,w) = V(y) + \frac{r}{2}\Phi(w)$. Moreover, if $\gamma(0) = 0$ then system (2.1), (2.2) with input $\delta \in \mathbb{R}^q$ and $\theta = 0$, $d = 0$ satisfies the ISS property,

(iii) that the subsystem (2.2) with inputs $(y,\delta) \in \mathbb{R}^n \times \mathbb{R}^q$ satisfies the p-ISS property; this is a consequence of (2.4). Moreover, if $\gamma(0) = 0$ then the subsystem (2.2) with inputs $(y,\delta) \in \mathbb{R}^n \times \mathbb{R}^q$ satisfies the ISS property.

However, due to the presence of $\theta \in \mathbb{R}^p$ and $d \in \mathbb{R}^m$ we cannot apply the feedback law $u = k(y)$ without some modification. Moreover, it should be noticed that since $\theta \in \mathbb{R}^p$ does not take values in a compact set, the "strength" of the interconnection between subsystems (2.1) and (2.2) is arbitrarily large and unknown.

Inequality (2.5) is automatically satisfied for appropriate $\Lambda > 0$, $\mu : \mathbb{R}^n \to (0,+\infty)$ when $\varphi(y,w) = \varphi(y,0)$ and $A(y,w) = A(y,0)$ for all $y \in \mathbb{R}^n$, $w \in \mathbb{R}^l$, i.e., when $\varphi$ and $A$ are independent of $w$. Moreover, inequality (2.6) is automatically satisfied for appropriate $\Lambda > 0$, $\mu : \mathbb{R}^n \to (0,+\infty)$.

Assumption (A) allows us to achieve robust stabilization of (2.1), (2.2) by means of an adaptive feedback law that depends only on $y \in \mathbb{R}^n$. Consider the DADS feedback law:

$$u = k(y) - C\mu^2(y)\frac{(b+\exp(z))^3}{a^3\beta^2}(\nabla V(y)g(y))^3$$
$$-C\left(|A(y,0)|^2 + |\varphi(y,0)|^2 + \mu(y) + 1\right)\frac{(b+\exp(z))^3}{a^3\beta^2}\nabla V(y)g(y) \quad (2.7)$$

$$\dot{z} = \Gamma\exp(-z)(V(y) - \varepsilon)^+ \quad , \quad z \in \mathbb{R} \quad (2.8)$$

where $\varepsilon, \Gamma > 0$, $a, \beta \in (0,1]$, $b, C \geq 1$ are parameters of the controller (constants) with $2a\beta < br$.

We call the controller (2.7), (2.8) a Deadzone-Adapted Disturbance Suppression (DADS) controller. The controller (2.7), (2.8) combines the use of deadzone (in (2.8)) and dynamic nonlinear damping (in (2.7)). The dynamic gain $z$ is being adapted by means of the update law (2.8). However, notice that $\dot{z}$ becomes zero, i.e., the adaptation stops, when the plant state $y$ enters the region defined by $V(y) \leq \varepsilon$. The deadzone prevents the dynamic gain $z$ from growing without bound in the case where a bounded disturbance is present.

The controller (2.7), (2.8) is simple: only one integrator is being used. The dynamic gain $z$ increases in order to overcome the effect of $\theta \in \mathbb{R}^p$ and $(d,\delta) \in \mathbb{R}^m \times \mathbb{R}^q$.

The following theorem clarifies the performance characteristics that the DADS controller (2.7), (2.8) can guarantee for the closed-loop system.

**Theorem 1:** *Suppose that Assumption (A) holds. Let $\varepsilon, \Gamma > 0$, $a, \beta \in (0,1]$, $b, C \geq 1$ with $2a\beta < br$ be given constants. Then there exist functions $\omega \in KL$, $\zeta \in K_\infty$, a non-decreasing function $\tilde{\zeta} \in C^0(\mathbb{R}_+; \mathbb{R}_+)$ and $B \in C^0(\mathbb{R}^{n+l+1} \times \mathbb{R}_+^3)$, such that for every $(y_0, w_0, z_0) \in \mathbb{R}^{n+l+1}$, $d \in L^\infty(\mathbb{R}_+; \mathbb{R}^m)$, $\theta \in L^\infty(\mathbb{R}_+; \mathbb{R}^p)$, $\delta \in L^\infty(\mathbb{R}_+; \mathbb{R}^q)$ the unique solution of the initial-value problem (2.1), (2.2), (2.7), (2.8) with initial condition $(y(0), w(0), z(0)) = (y_0, w_0, z_0)$ satisfies the following estimates for all $t \geq 0$:*

$$\limsup_{t \to +\infty}(V(y(t))) \leq \varepsilon, \tag{2.9}$$

$$\limsup_{t \to +\infty}(\Phi(w(t))) \leq \tilde{\zeta}\left(\varepsilon + \limsup_{t \to +\infty}(|\delta(t)|)\right), \tag{2.10}$$

$$U(y(t), w(t)) \leq \omega(U(y_0, w_0), t) + \zeta\left(\chi(\|d\|_\infty, \|\delta\|_\infty, \|\theta\|_\infty, \exp(z_0))\right), \tag{2.11}$$

$$z_0 \leq z(t) \leq \lim_{s \to +\infty}(z(s)) \leq B(y_0, w_0, z_0, \|d\|_\infty, \|\delta\|_\infty, \|\theta\|_\infty). \tag{2.12}$$

*where*

$$\chi(s_1, s_2, s_3, s_4) := \frac{r}{2}\gamma(s_2) + a \frac{s_1^2 + s_1^4 + \left((s_3 - b - s_4)^+\right)^2 + \left((s_3 - b - s_4)^+\right)^4 + 2\beta\Lambda}{b + s_4} \tag{2.13}$$

*for all $s_1, s_2, s_3, s_4 \geq 0$ and*

$$U(y, w) := V(y) + \frac{r}{2}\Phi(w) \tag{2.14}$$

*for all $(y, w) \in \mathbb{R}^{n+l}$. Moreover, if $\gamma(0) = 0$ then $\tilde{\zeta}(0) = 0$. Finally, if $\gamma(0) = \Lambda = 0$ and $\lim_{t \to +\infty}(d(t)) = 0$, $\lim_{t \to +\infty}(\delta(t)) = 0$ and $\theta$ is constant then one of the following holds: either $\lim_{t \to +\infty}(|y(t)|) = \lim_{t \to +\infty}(|w(t)|) = 0$ or $|\theta| > b$ and $\lim_{s \to +\infty}(z(s)) < \ln(|\theta| - b)$.*

Estimate (2.11) shows the p-IOS property for system (2.1), (2.2), (2.7), (2.8) with constant $\theta$, $(d, \delta) \in \mathbb{R}^m \times \mathbb{R}^q$ as input and $(y, w) \in \mathbb{R}^n \times \mathbb{R}^l$ as output. Estimates (2.11) and (2.12) show the p-UBIBS property for system (2.1), (2.2), (2.7), (2.8) with $(d, \delta, \theta) \in \mathbb{R}^m \times \mathbb{R}^q \times \mathbb{R}^p$ as input.

It can be shown (see the proof of Theorem 2 below) that Assumption (A) is automatically satisfied when the following assumptions are valid:

**(B1)** *Subsystem (2.2) with inputs $(y, \delta) \in \mathbb{R}^n \times \mathbb{R}^q$ satisfies the ISS property,*

**(B2)** *There exists a smooth mapping $k: \mathbb{R}^n \to \mathbb{R}$ with $k(0) = 0$ such that $0 \in \mathbb{R}^n$ is globally asymptotically stable for the closed-loop system (2.1) with $\theta = 0$, $d = 0$ and $u = k(y)$.*

On the other hand, it should be noticed that assumptions (B1), (B2) are not equivalent with assumption (A). For example, the system $\dot{y} = u + \theta w^2$, $\dot{w} = -(w^2 - 1 - y^2)w$ with $y, u, w, \theta \in \mathbb{R}$ satisfies assumption (A) with $k(y) = -y - y^3$, $V(y) = y^2/2$, $Q(y) = y^2 + y^4$, $r = 1$, $\Lambda = 0$, $\mu(y) \equiv 2$, $\Phi(w) = w^2/2$, $R(w) = w^4/2$ and $\gamma(s) \equiv 1$. However, assumption (B1) is not valid since the subsystem $\dot{w} = -(w^2 - 1 - y^2)w$ with $y \in \mathbb{R}$ as input is not ISS. To see this, notice that $0 \in \mathbb{R}$ is not globally asymptotically stable when $y \equiv 0$ (because there are the non-zero equilibrium points $w = \pm 1$).

Exploiting assumptions (B1), (B2), we obtain the following result.

**Theorem 2:** *Suppose that assumptions (B1), (B2) are valid for system (2.1), (2.2). Then there exist smooth mappings $V : \mathbb{R}^n \to \mathbb{R}_+$, $\Phi : \mathbb{R}^l \to \mathbb{R}_+$, $\mu : \mathbb{R}^n \to (0, +\infty)$, $V, \Phi$ being positive definite and radially unbounded, constants $r, \Lambda > 0$ and a non-decreasing function $\gamma \in C^0(\mathbb{R}_+; \mathbb{R}_+)$ such that for every $\varepsilon, \Gamma > 0$, $a, \beta \in (0,1]$, $b, C \geq 1$ with $2a\beta < br$ there exist functions $\omega \in KL$, $\zeta \in K_\infty$, a non-decreasing function $\tilde{\zeta} \in C^0(\mathbb{R}_+; \mathbb{R}_+)$ and $B \in C^0(\mathbb{R}^{n+l+1} \times \mathbb{R}_+^3)$, with the property that for every $(y_0, w_0, z_0) \in \mathbb{R}^{n+l+1}$, $d \in L^\infty(\mathbb{R}_+; \mathbb{R}^m)$, $\theta \in L^\infty(\mathbb{R}_+; \mathbb{R}^p)$, $\delta \in L^\infty(\mathbb{R}_+; \mathbb{R}^q)$ the unique solution of the initial-value problem (2.1), (2.2), (2.7), (2.8) with initial condition $(y(0), w(0), z(0)) = (y_0, w_0, z_0)$ satisfies estimates (2.9), (2.10), (2.11), (2.12) for all $t \geq 0$.*

For the proof of Theorem 1, we use the following technical lemma, proven in Section 4.

**Lemma 1:** *Let $\rho \in K_\infty$, $T, \varepsilon > 0$ be given. Then there exists a positive non-increasing function $c_\varepsilon : \mathbb{R}_+ \to (0, +\infty)$ such that for every absolutely continuous function $V : [0, T) \to \mathbb{R}_+$ and for every non-increasing function $\alpha \in C^0([0, T); \mathbb{R}_+)$ for which the following differential inequality holds for $t \in [0, T)$ a.e.*

$$\dot{V}(t) \leq -\rho(V(t)) + \alpha(t) \tag{2.15}$$

*the following estimate holds for all $t_0 \in [0, T)$ and $t \in [t_0, T)$*

$$V(t) \leq \min\left(s, s \exp(-c_\varepsilon(s)(t - t_0)) + \frac{\varepsilon}{2} + \frac{\alpha(t_0)}{c_\varepsilon(s)}\right) \tag{2.16}$$

*where $s = V(0) + \rho^{-1}(\alpha(0))$.*

# 3. DADS is not limited to finite-dimensional systems

The partial-state DADS controller may be applied even in cases where the unmeasured state component belongs to an infinite-dimensional space. To illustrate this point, consider the following system

$$w_t(x) = pw_{xx}(x) + \theta_1 K(x, y)$$
$$w(0) = w(1) = 0 \quad (3.1)$$
$$\dot{y} = u + \theta_2 L(w) + d$$

where $x \in [0,1]$ is the spatial variable, $y \in \mathbb{R}$, $w \in H^2(0,1) \cap H_0^1(0,1)$ are the state components, $u \in \mathbb{R}$ is the control input, $d \in \mathbb{R}$ is the disturbance, $\theta = (\theta_1, \theta_2) \in \mathbb{R}^2$ are unknown possibly time-varying parameters, $p > 0$ is an unknown constant parameter, $K: [0,1] \times \mathbb{R} \to \mathbb{R}$ is an unknown function and $L: H^2(0,1) \cap H_0^1(0,1) \to \mathbb{R}$ is an unknown functional that satisfy the following inequalities for all $x \in [0,1]$, $y \in \mathbb{R}$ and $w \in H^2(0,1) \cap H_0^1(0,1)$

$$|L(w)| \leq \|w\|, \quad |K(x, y)| \leq |y| \quad (3.2)$$

where $\|w\|$ is the standard norm in $L^2(0,1)$. Of course, it should be noted that in order to be able to guarantee existence and uniqueness of solutions of an initial-boundary value problem for system (3.1) additional assumptions are required for the function $K: [0,1] \times \mathbb{R} \to \mathbb{R}$ and the functional $L: H^2(0,1) \cap H_0^1(0,1) \to \mathbb{R}$. However, next we do not study existence/uniqueness issues for system (3.1) and we focus only on the control issues.

We seek a robust adaptive controller that depends only on $y \in \mathbb{R}$ and guarantees robust practical regulation of the state of system (3.1).

The feedback stabilization problem described above is highly non-trivial. First of all, it should be noticed that almost nothing is known about the dynamics of the unmeasured state component $w$: $\theta_1, \theta_2 \in \mathbb{R}$, $p > 0$ are unknown constant parameters, $K: [0,1] \times \mathbb{R} \to \mathbb{R}$ and $L: H^2(0,1) \cap H_0^1(0,1) \to \mathbb{R}$ are unknown functionals that satisfy (3.2). Secondly, it should be emphasized that even when $d = 0$ and $\theta = (\theta_1, \theta_2) \in \mathbb{R}^2$ is constant, the equilibrium point $(y, w) = (0, 0)$ of the open-loop system (3.1) can be unstable and system (3.1) can have exponentially increasing (in norm) solutions. For example, when $L(w) = -\int_0^1 w(x)dx$, $K(x, y) = \frac{x^2 - x - 2p}{1 + 2p} y$, $6(1 + 2p) = \theta_1 \theta_2$ and $d = 0$ the reader can verify that (3.2) is valid and $y(t) = \frac{\theta_2}{6} \exp(t), w(t, x) = \exp(t)x(x - 1)$ is an exponentially increasing (in norm) solution of the open-loop system (3.1).

System (3.1) is not a straightforward extension of the finite-dimensional case (2.1), (2.2). Indeed, the reader can notice that (3.1) allows the unmeasured state $w$ to be directly affected by the unknown possibly time-varying parameters $\theta = (\theta_1, \theta_2) \in \mathbb{R}^2$. Moreover, it should be noted that the dynamics of the unmeasured state contain an additional unknown parameter, namely the diffusion coefficient $p > 0$, which is not denoted in the same way as the other unknown parameters (e.g., we do not use the notation $\theta_3$ for the diffusion coefficient) because it is a constant parameter. Therefore, the DADS feedback design procedure that is described below is slightly different from the design described in the previous section. Moreover, the obtained results are also slightly different: we obtain stability estimates with different powers from the powers appearing in (2.13). Consequently, the present application shows that partial state DADS feedback may be possible in cases more demanding than (2.1), (2.2).

In order to perform the DADS feedback design, we use the functionals

$$\Phi(w) = \|w\|^2 / 2 \tag{3.3}$$

$$V(y) = y^2 / 2 \tag{3.4}$$

Using (3.1), (3.2), (3.3), integration by parts and Wirtinger's inequality ($\|w'\|^2 \geq \pi^2 \|w\|^2$ for $w \in H^2(0,1) \cap H_0^1(0,1)$) we get for all $w \in H^2(0,1) \cap H_0^1(0,1)$, $y \in \mathbb{R}$:

$$\dot{\Phi} = -p\|w_x\|^2 + \theta_1 \int_0^1 K(x,y) w(x) dx \leq -p\pi^2 \|w\|^2 + |\theta_1| |y| \|w\| \tag{3.5}$$

Using the fact $|\theta_1| \leq |\theta|$ and the inequalities

$$|\theta| |y| \|w\| \leq \frac{p\pi^2}{4} \|w\|^2 + \frac{|\theta|^2}{p\pi^2} y^2, \quad \frac{|\theta|^2}{p\pi^2} y^2 \leq \frac{|\theta|^4}{4 p^2 \pi^4} + y^4$$

we get from (3.3), (3.5) for all $w \in H^2(0,1) \cap H_0^1(0,1)$, $y \in \mathbb{R}$:

$$\dot{\Phi} \leq -\frac{3 p \pi^2}{2} \Phi(w) + \frac{|\theta|^4}{4 p^2 \pi^4} + y^4 \tag{3.6}$$

Using the inequality $|\theta_1| |y| \|w\| \leq \frac{p\pi^2}{2} \|w\|^2 + \frac{|\theta_1|^2}{2 p \pi^2} y^2$, we also get from (3.3), (3.4), (3.5) for all $w \in H^2(0,1) \cap H_0^1(0,1)$, $y \in \mathbb{R}$:

$$\dot{\Phi} \leq -p\pi^2 \Phi(w) + \frac{|\theta_1|^2}{p\pi^2} V(y) \tag{3.7}$$

Let $b \geq 1 \geq a > 0$ be given (arbitrary) constants. Using (3.1), (3.2), (3.4) and the facts $|\theta_2| \leq |\theta|$, $|\theta| \leq (|\theta| - b - \exp(z))^+ + (b + \exp(z))$, we get for all $w \in H^2(0,1) \cap H_0^1(0,1)$, $y \in \mathbb{R}$, $z \in \mathbb{R}$, $d \in \mathbb{R}$:

$$\begin{aligned}\dot{V} &= yu + \theta_2 y L(w) + yd \\ &\leq yu + |\theta| |y| \|w\| + yd \\ &\leq yu + (|\theta| - b - \exp(z))^+ |y| \|w\| + (b + \exp(z)) |y| \|w\| + yd\end{aligned} \tag{3.8}$$

Using the inequalities

$$yd \leq \frac{ad^2}{b+\exp(z)} + \frac{b+\exp(z)}{4a}y^2, \quad \frac{(b+\exp(z))^3}{ap\pi^2}y^2 \leq \frac{a}{p^2(b+\exp(z))} + \frac{(b+\exp(z))^7}{4a^3\pi^4}y^4$$

$$(b+\exp(z))|y|\|w\| \leq \frac{ap\pi^2}{4(b+\exp(z))}\|w\|^2 + \frac{(b+\exp(z))^3}{ap\pi^2}y^2$$

$$\left(|\theta|-b-\exp(z)\right)^+|y|\|w\| \leq \frac{ap\pi^2}{4(b+\exp(z))}\|w\|^2 + \frac{(b+\exp(z))}{ap\pi^2}\left(\left(|\theta|-b-\exp(z)\right)^+\right)^2 y^2$$

$$\frac{(b+\exp(z))}{ap\pi^2}\left(\left(|\theta|-b-\exp(z)\right)^+\right)^2 y^2 \leq \frac{a}{p^2(b+\exp(z))} + \frac{(b+\exp(z))^3}{4a^3\pi^4}\left(\left(|\theta|-b-\exp(z)\right)^+\right)^4 y^4$$

$$\frac{(b+\exp(z))^3}{4a^3\pi^4}\left(\left(|\theta|-b-\exp(z)\right)^+\right)^4 y^4 \leq a\frac{\left(\left(|\theta|-b-\exp(z)\right)^+\right)^8}{b+\exp(z)} + \frac{(b+\exp(z))^7}{64a^6\pi^8}y^8$$

we obtain from (3.3), (3.8) for all $w \in H^2(0,1) \cap H_0^1(0,1)$, $y \in \mathbb{R}$, $z \in \mathbb{R}$, $d \in \mathbb{R}$:

$$\dot{V} \leq yu + \frac{(b+\exp(z))^7}{64a^6\pi^8}y^8 + \frac{(b+\exp(z))^7}{4a^3\pi^4}y^4 + \frac{b+\exp(z)}{4a}y^2 \\ + a\frac{d^2 + \left(\left(|\theta|-b-\exp(z)\right)^+\right)^8 + 2p^{-2} + p\pi^2\Phi(w)}{b+\exp(z)} \tag{3.9}$$

The fact that $b \geq 1 \geq a > 0$ in conjunction with (3.9) implies that the following inequality holds for all $w \in H^2(0,1) \cap H_0^1(0,1)$, $y \in \mathbb{R}$, $z \in \mathbb{R}$, $d \in \mathbb{R}$:

$$\dot{V} \leq yu + \frac{(b+\exp(z))^7}{4a^6}\left(y^8 + y^4 + y^2\right) + a\frac{d^2 + \left(\left(|\theta|-b-\exp(z)\right)^+\right)^8 + 2p^{-2} + p\pi^2\Phi(w)}{b+\exp(z)} \tag{3.10}$$

Let $c \geq 1$ be a given (arbitrary) constant. By virtue of (3.4), (3.10) the feedback law

$$u = -c\left(1 + \frac{(b+\exp(z))^7}{4a^6}\right)\left(y^7 + y^3 + y\right) \tag{3.11}$$

guarantees the following differential inequality for all $w \in H^2(0,1) \cap H_0^1(0,1)$, $y \in \mathbb{R}$, $z \in \mathbb{R}$, $d \in \mathbb{R}$:

$$\dot{V} \leq -2cV(y) - cy^4 - cy^8 + a\frac{d^2 + \left(\left(|\theta|-b-\exp(z)\right)^+\right)^8 + 2p^{-2} + p\pi^2\Phi(w)}{b+\exp(z)} \tag{3.12}$$

The sole reason for the appearance of the seventh powers in (3.11) and the eight powers in (3.9), whereas such powers do not appear in (2.7), is that the $w-$system (3.1) includes a dependency on the unknown parameter $\theta$, which has to be dominated in the DADS design, while the $w-$system (2.2) has no such $\theta-$dependency.

Let $\varepsilon, \Gamma > 0$ be given (arbitrary) constants. The DADS feedback law (3.11) with (2.8) can guarantee the practical regulation of the plant state and boundedness of solutions for every bounded parameter $\theta = (\theta_1, \theta_2) \in \mathbb{R}^2$, for every bounded disturbance $d \in \mathbb{R}$ and for every initial condition. In order to understand how all these properties are achieved, we provide the following result.

**Theorem 3:** *Let $\varepsilon, \Gamma > 0$, $c \geq 1$, $b \geq 1 \geq a > 0$ be given constants. Then every solution*

$$w \in C^0(\mathbb{R}_+; L^2(0,1)) \cap C^1((0,+\infty); L^2(0,1)) \cap C^0((0,+\infty); H^2(0,1) \cap H_0^1(0,1)),$$

$z \in C^1(\mathbb{R}_+; \mathbb{R})$, $y \in C^0(\mathbb{R}_+; \mathbb{R})$ *being absolutely continuous on every compact interval of $\mathbb{R}_+$*

*of the closed-loop system (3.1) with (3.11), (2.8) which satisfies $w_t(x) = p w_{xx}(x) + \theta_1 K(x, y)$, (3.10), (2.8) for all $t \in (0, +\infty)$ and $\dot{y} = u + \theta_2 L(w) + d$ for $t \in (0, +\infty)$ a.e. for some $d \in L^\infty(\mathbb{R}_+)$, $\theta \in C^0(\mathbb{R}_+; \mathbb{R}^2) \cap L^\infty(\mathbb{R}_+; \mathbb{R}^2)$ satisfies the following estimates for all $t \geq 0$*

$$\|w[t]\|^2 + y^2(t) \leq \exp(-\kappa t)(\|w[0]\|^2 + y^2(0)) + 2a \frac{\|d\|_\infty^2 + \|\theta\|_\infty^8 + 2p^{-2}}{b\kappa} + \frac{\|\theta\|_\infty^4}{2p^2\pi^4\kappa} \quad (3.13)$$

$$z(0) \leq z(t) \leq \ln\left(\exp(z(0)) + \bar{B}\left(\|d\|_\infty^2 + \|\theta\|_\infty^4 + \|\theta\|_\infty^8 + 2p^{-2} + \|w[0]\|^2 + y^2(0)\right)\right) \quad (3.14)$$

$$\limsup_{t \to +\infty}(\|w[t]\|) \leq \frac{\sqrt{2\varepsilon}}{p\pi^2} \limsup_{t \to +\infty}(|\theta_1(t)|) \quad (3.15)$$

$$\limsup_{t \to +\infty}(|y(t)|) \leq \sqrt{2\varepsilon} \quad (3.16)$$

*where* $\kappa = \min\left(\frac{p\pi^2}{2}, 2c\right)$, $\bar{K} = \frac{4p\pi^2(b\kappa + ap\pi^2) + 2b\kappa p^2\pi^4 + b}{4p\pi^2 b\kappa}$ *and* $\bar{B} = \frac{a\bar{K}(\varepsilon\Gamma + 2c) + \varepsilon c\Gamma}{4c^2\varepsilon}$.

It should be noticed that estimate (3.13) shows an infinite-dimensional version of the p-IOS property of the closed-loop system (3.1) with (3.11), (2.8) and output $(w, y) \in L^2(0,1) \times \mathbb{R}$. Moreover, the combination of (3.13) and (3.14) shows an infinite-dimensional version of the p-UBIBS property. On the other hand, when $\theta$ is constant, the asymptotic estimates (3.15) and (3.16) show the zero p-OAG for the closed-loop system (3.1) with (3.11), (2.8) and output $(w, y) \in L^2(0,1) \times \mathbb{R}$. Estimate (3.16) alone shows the zero p-OAG property for the closed-loop system (3.1) with (3.11), (2.8) and output $y \in \mathbb{R}$ even in the case of time-varying $\theta$.

Again, it should be noted that the difference between (3.13) and (2.11) (different powers in $\|d\|_\infty$ and $\|\theta\|_\infty$) is explained by the fact that (3.1) allows the unmeasured state $w$ to be directly affected by the unknown parameters, while this is not the case in (2.1), (2.2). Due to this structural difference, the DADS feedback design procedure is slightly different here from the design described in Section 2 and the obtained results are slightly different. It should be noticed that the infinite-dimensional character of the unmeasured state component $w$ plays no role in this. To understand this point, we note that for every $\varepsilon, \Gamma, p > 0$, $c \geq 1$, $b \geq 1 \geq a > 0$ exactly the same estimates (3.13), (3.14), (3.15), (3.16) - with $\|w\|$ replaced by $|w|$ - are obtained for the finite-dimensional closed-loop system

$$\dot{w} = -p\pi^2 w + \theta_1 y$$
$$\dot{y} = u + \theta_2 w + d \tag{3.17}$$
$$y, w, u, d \in \mathbb{R}, \theta = (\theta_1, \theta_2) \in \mathbb{R}^2$$

with (3.10), (2.8). System (3.17) is not a system of the form (2.1), (2.2) since $\theta$ appears in the equation for $\dot{w}$. That's why we noted above that the present application shows that partial state DADS feedback may be possible in cases more demanding than (2.1), (2.2).

The proof of Theorem 3 is similar to the proof of Theorem 1. We next give a brief sketch of the proof of Theorem 3.

**Sketch of Proof of Theorem 3:** We repeat the steps of Theorem 1.

1) By utilizing the functional
$$U(y, w) = \left( \|w\|^2 + y^2 \right) / 2 \tag{3.18}$$

defined for all $w \in L^2(0,1)$, $y \in \mathbb{R}$, and exploiting (3.6), (3.12) and the facts that $b \geq 1 \geq a > 0$, $c \geq 1$, we obtain for all $w \in H^2(0,1) \cap H_0^1(0,1)$, $y \in \mathbb{R}$, $z \in \mathbb{R}$, $d \in \mathbb{R}$:

$$\begin{aligned}\dot{U} &\leq -\frac{p\pi^2}{2}\Phi(w) - 2cV(y) + b^{-1}a\left(d^2 + \left((|\theta|-b)^+\right)^8 + 2p^{-2}\right) + \frac{|\theta|^4}{4p^2\pi^4} \\ &\leq -\min\left(\frac{p\pi^2}{2}, 2c\right)U(y,w) + b^{-1}a\left(d^2 + \left((|\theta|-b)^+\right)^8 + 2p^{-2}\right) + \frac{|\theta|^4}{4p^2\pi^4}\end{aligned} \tag{3.19}$$

Differential inequality (3.19) and definition (3.18) show that every solution $w \in C^0(\mathbb{R}_+; L^2(0,1)) \cap C^1((0,+\infty); L^2(0,1)) \cap C^0((0,+\infty); H^2(0,1) \cap H_0^1(0,1))$, $z \in C^1(\mathbb{R}_+; \mathbb{R})$, $y \in C^0(\mathbb{R}_+; \mathbb{R})$ being absolutely continuous of the closed-loop system (3.1) with (3.11), (2.8) which satisfies $w_t(x) = pw_{xx}(x) + \theta_1 K(x, y)$, (3.10), (2.8) for all $t \in (0, +\infty)$ and $\dot{y} = u + \theta_2 L(w) + d$ for $t \in (0, +\infty)$ a.e. for some $d \in L^\infty(\mathbb{R}_+)$, $\theta \in C^0(\mathbb{R}_+; \mathbb{R}^2) \cap L^\infty(\mathbb{R}_+; \mathbb{R}^2)$ satisfies estimate (3.13) for all $t \geq 0$.

2) Next we show that $z(t)$ is bounded. Exploiting inequality (3.12), definition (3.3) and estimate (3.13) we get for $t \geq 0$ a.e.:

$$\frac{d}{dt}(V(y(t))) \leq -2cV(y(t)) + R(b + \exp(z(t)))^{-1} \tag{3.20}$$

where $R = a\bar{K}\left(\|d\|_\infty^2 + \|\theta\|_\infty^4 + \|\theta\|_\infty^8 + 2p^{-2} + \|w[0]\|^2 + y^2(0)\right)$. Therefore, by using Lemma 4.1 in [11] we can guarantee estimate (3.14) for all $t \geq 0$.

3) We next show estimate (2.9). Estimate (3.14) and the fact that $z(t)$ is non-decreasing guarantees that the function $z(t)$ has a finite limit as $t \to +\infty$. This implies that the function $\exp(z(t))$ has a finite limit as $t \to +\infty$. Moreover, the facts that $d \in L^\infty(\mathbb{R}_+)$, $\theta \in L^\infty(\mathbb{R}_+; \mathbb{R}^2)$, $\|w[\cdot]\|, y, u \in L^\infty(\mathbb{R}_+)$ and (3.1), (3.11) imply that $\frac{d}{dt}(V(y(t)))$ is of class $L^\infty(\mathbb{R}_+)$. It follows that

the function $\frac{d}{dt}(\exp(z(t))) = \Gamma(V(y(t)) - \varepsilon)^+$ is uniformly continuous, in addition to $\int_0^{+\infty} \frac{d}{dt}(\exp(z(t)))\,dt = \lim_{t \to +\infty}(\exp(z(t))) - \exp(z(0)) < +\infty$. From Barbălat's Lemma (see [12, 7]), we have:

$$\lim_{t \to +\infty}\left(\frac{d}{dt}(\exp(z(t)))\right) = \lim_{t \to +\infty}\left(\Gamma(V(y(t)) - \varepsilon)^+\right) = 0 \tag{3.21}$$

Therefore, estimate (2.9) holds.

4) We finally provide an asymptotic estimate for the unmeasured state $w$. Differential inequality (3.7) in conjunction with (2.9) show that

$$\limsup_{t \to +\infty}\left(\Phi(w[t])\right) \le \frac{\varepsilon}{p^2 \pi^4} \limsup_{t \to +\infty}\left(|\theta_1(t)|^2\right)$$

The above estimate in conjunction with definitions (3.3), (3.4) and estimate (2.9) show that the asymptotic estimates (3.15), (3.16) are valid. The proof is complete. ◁

## 4. Proofs of Main Results

We first provide the proof of Lemma 1.

**Proof of Lemma 1:** We first construct the positive, non-increasing function $c_\varepsilon : \mathbb{R}_+ \to (0, +\infty)$ so that

$$c_\varepsilon(\tau)\left(\tau - \frac{\varepsilon}{2}\right)^+ \le \rho(\tau) \text{ for all } \tau \ge 0 \tag{4.1}$$

In fact, we can give a formula for $c_\varepsilon : \mathbb{R}_+ \to (0, +\infty)$: the reader can verify that the function defined by the formula

$$c_\varepsilon(\tau) := \begin{cases} \min\left(1, \inf\left\{\frac{2\rho(l)}{2l - \varepsilon} : \frac{\varepsilon}{2} < l \le \tau\right\}\right) & \text{for } \tau > \varepsilon/2 \\ 1 & \text{for } \tau \in [0, \varepsilon/2] \end{cases} \tag{4.2}$$

is a positive, non-increasing function $c_\varepsilon : \mathbb{R}_+ \to (0, +\infty)$ that satisfies (4.1).

Define the function for $t \in [0, T)$

$$\lambda(t) = \frac{1}{2}\left(\left(V(t) - \rho^{-1}(\alpha(0))\right)^+\right)^2 \tag{4.3}$$

Then we get from (2.15), (4.3) the following differential inequality for $t \in [0, T)$ a.e.:

$$\dot{\lambda}(t) = \left(V(t) - \rho^{-1}(\alpha(0))\right)^+ \dot{V}(t) \leq \left(V(t) - \rho^{-1}(\alpha(0))\right)^+ \left(-\rho(V(t)) + \alpha(t)\right)$$
$$\leq -\left(V(t) - \rho^{-1}(\alpha(0))\right)^+ \left(\rho(V(t)) - \alpha(0)\right) \leq 0$$

Therefore, the function $\lambda$ defined by (4.3) is non-increasing on $[0,T)$ and we obtain the estimate for $t \in [0,T)$ (recall that $s = V(0) + \rho^{-1}(\alpha(0))$):

$$V(t) \leq s \tag{4.4}$$

It follows from (2.15) and (4.1) that the following differential inequality holds for $t \in [0,T)$ a.e.

$$\dot{V}(t) \leq -c_\varepsilon(V(t))\left(V(t) - \frac{\varepsilon}{2}\right)^+ + \alpha(t) \tag{4.5}$$

Since $c_\varepsilon : \mathbb{R}_+ \to (0,+\infty)$ is non-increasing we get from (4.4) and (4.5) that the following differential inequality holds for $t \in [0,T)$ a.e.

$$\dot{V}(t) \leq -c_\varepsilon(s)\left(V(t) - \frac{\varepsilon}{2}\right)^+ + \alpha(t) \tag{4.6}$$

Consequently, we conclude from (4.6) and the fact that $\alpha \in C^0([0,T); \mathbb{R}_+)$ is non-increasing that the following differential inequality holds for all $t_0 \in [0,T)$ and $t \in [t_0, T)$ a.e.

$$\dot{V}(t) \leq -c_\varepsilon(s)\left(V(t) - \frac{\varepsilon}{2}\right) + \alpha(t_0) \tag{4.7}$$

Integrating directly (using the integrating factor $\exp(c_\varepsilon(s)t)$) the differential inequality (4.7), we obtain for all $t_0 \in [0,T)$ and $t \in [t_0,T)$:

$$V(t) \leq \exp(-c_\varepsilon(s)(t-t_0))V(t_0) + \left(\frac{\varepsilon}{2} + \frac{\alpha(t_0)}{c_\varepsilon(s)}\right)\left(1 - \exp(-c_\varepsilon(s)(t-t_0))\right) \tag{4.8}$$

Combining (4.8) and (4.4) we get (2.16) for all $t \in [t_0, T)$. The proof is complete. ◁

We can now give the proof of Theorem 1.

**Proof of Theorem 1:** First we show some useful inequalities.

We have for all $(y, w, \theta) \in \mathbb{R}^{n+l+p}$:

$$\nabla V(y)g(y)\varphi'(y,w)\theta \leq |\nabla V(y)g(y)||\varphi(y,w)||\theta|$$
$$\leq |\nabla V(y)g(y)||\varphi(y,w)-\varphi(y,0)||\theta| + |\nabla V(y)g(y)||\varphi(y,0)||\theta|$$
$$\leq \frac{b+\exp(z)}{4a\beta}\mu(y)(\nabla V(y)g(y))^2|\theta|^2 + \frac{a\beta|\varphi(y,w)-\varphi(y,0)|^2}{(b+\exp(z))\mu(y)}$$
$$+ |\nabla V(y)g(y)||\varphi(y,0)||\theta|$$

The inequality $|\theta| \leq (|\theta|-b-\exp(z))^+ + b+\exp(z)$ implies that $|\theta|^2 \leq 2((|\theta|-b-\exp(z))^+)^2 + 2(b+\exp(z))^2$. Thus, we obtain from the above inequality for all $(y,w,\theta) \in \mathbb{R}^{n+l+p}$:

$$\nabla V(y)g(y)\varphi'(y,w)\theta \leq \frac{b+\exp(z)}{2a\beta}\mu(y)(\nabla V(y)g(y))^2 \left((|\theta|-b-\exp(z))^+\right)^2$$
$$+ \frac{(b+\exp(z))^3}{2a\beta}\mu(y)(\nabla V(y)g(y))^2 + \frac{a\beta|\varphi(y,w)-\varphi(y,0)|^2}{(b+\exp(z))\mu(y)}$$
$$+ |\nabla V(y)g(y)||\varphi(y,0)|(|\theta|-b-\exp(z))^+$$
$$+ |\nabla V(y)g(y)||\varphi(y,0)|(b+\exp(z))$$

Completing the squares, we get for all $(y,w,\theta) \in \mathbb{R}^{n+l+p}$:

$$\nabla V(y)g(y)\varphi'(y,w)\theta \leq a \frac{\left((|\theta|-b-\exp(z))^+\right)^4 + \left((|\theta|-b-\exp(z))^+\right)^2}{b+\exp(z)}$$
$$+ \frac{a\beta\left(|\varphi(y,w)-\varphi(y,0)|^2 + |\varphi(y,0)|^2\right)}{(b+\exp(z))\mu(y)} + \frac{(b+\exp(z))}{4a}(\nabla V(y)g(y))^2|\varphi(y,0)|^2 \quad (4.9)$$
$$+ \left(\frac{\mu(y)}{4a^2\beta}(\nabla V(y)g(y))^2 + 3\right)\frac{(b+\exp(z))^3}{4a\beta}\mu(y)(\nabla V(y)g(y))^2$$

We also have for all $(y,w,d) \in \mathbb{R}^{n+l+m}$:

$$\nabla V(y)g(y)A'(y,w)d = \nabla V(y)g(y)A'(y,0)d + \nabla V(y)g(y)(A(y,w)-A(y,0))'d$$
$$\leq |\nabla V(y)g(y)||A(y,0)||d| + |\nabla V(y)g(y)||A(y,w)-A(y,0)||d|$$

Using the facts that

$$|\nabla V(y)g(y)||A(y,0)||d| \leq \frac{b+\exp(z)}{4a}|A(y,0)|^2(\nabla V(y)g(y))^2 + \frac{a|d|^2}{b+\exp(z)}$$

$$|\nabla V(y)g(y)||A(y,w)-A(y,0)||d| \le \frac{a\beta|A(y,w)-A(y,0)|^2}{(b+\exp(z))\mu(y)}$$

$$+\frac{(b+\exp(z))\mu(y)}{4a\beta}(\nabla V(y)g(y))^2|d|^2$$

we get for all $(y,w,d) \in \mathbb{R}^{n+l+m}$:

$$\nabla V(y)g(y)A'(y,w)d \le \frac{(b+\exp(z))\mu(y)}{4a\beta}(\nabla V(y)g(y))^2|d|^2$$

$$+\frac{b+\exp(z)}{4a}|A(y,0)|^2(\nabla V(y)g(y))^2 + \frac{a|d|^2}{b+\exp(z)} + \frac{a\beta|A(y,w)-A(y,0)|^2}{(b+\exp(z))\mu(y)}$$

Completing the squares, we obtain for all $(y,w,d) \in \mathbb{R}^{n+l+m}$:

$$\nabla V(y)g(y)A'(y,w)d \le \frac{a|d|^2+a|d|^4}{b+\exp(z)} + \frac{a\beta|A(y,w)-A(y,0)|^2}{(b+\exp(z))\mu(y)}$$

$$+\left(\frac{(b+\exp(z))^2\mu^2(y)}{16a^2\beta^2}(\nabla V(y)g(y))^2 + |A(y,0)|^2\right)\frac{b+\exp(z)}{4a}(\nabla V(y)g(y))^2 \quad (4.10)$$

Combining (4.9), (4.10) and (2.5), (2.6), we obtain:

$$\nabla V(y)g(y)\left(\varphi'(y,w)\theta + A'(y,w)d\right) \le \left(|A(y,0)|^2 + |\varphi(y,0)|^2\right)\frac{b+\exp(z)}{4a}(\nabla V(y)g(y))^2$$

$$+a\frac{|d|^2+|d|^4+\left((|\theta|-b-\exp(z))^+\right)^2+\left((|\theta|-b-\exp(z))^+\right)^4}{b+\exp(z)}$$

$$+\frac{a\beta(R(w)+Q(y)+2\Lambda)}{(b+\exp(z))} + \frac{a\beta}{(b+\exp(z))}(\nabla V(y)g(y))^2$$

$$+\left(\frac{5\mu(y)}{16a^2\beta}(\nabla V(y)g(y))^2 + 3\right)\frac{(b+\exp(z))^3}{4a\beta}\mu(y)(\nabla V(y)g(y))^2$$

Since $b+\exp(z) \ge b$ and $(b+\exp(z))^3 \ge b+\exp(z)$ (a consequence of the fact that $b \ge 1$), $a, \beta \in (0,1]$, we get for all $(y,w,d,\theta) \in \mathbb{R}^{n+l+m+p}$:

$$\nabla V(y)g(y)\left(\varphi'(y,w)\theta + A'(y,w)d\right) \le \frac{a\beta(R(w)+Q(y)+2\Lambda)}{b+\exp(z)}$$

$$+a\frac{|d|^2+|d|^4+\left((|\theta|-b-\exp(z))^+\right)^2+\left((|\theta|-b-\exp(z))^+\right)^4}{b+\exp(z)} \quad (4.11)$$

$$+\left(|A(y,0)|^2 + |\varphi(y,0)|^2 + \mu(y) + 1 + \mu^2(y)(\nabla V(y)g(y))^2\right)\frac{(b+\exp(z))^3}{a^3\beta^2}(\nabla V(y)g(y))^2$$

Using definition (2.14), (2.1), (2.2) and (2.7) we get for all $(y,w,d,\theta,\delta) \in \mathbb{R}^{n+l+m+p+q}$:

$$\dot{U} = \nabla V(y)(f(y)+g(y)k(y)) + \frac{r}{2}\nabla\Phi(w)h(y,w,\delta)$$
$$+\nabla V(y)g(y)(\varphi'(y,w)\theta + A'(y,w)d) - C\frac{(b+\exp(z))^3}{a^3\beta^2}\mu^2(y)(\nabla V(y)g(y))^4 \quad (4.12)$$
$$-C\left(|A(y,0)|^2 + |\varphi(y,0)|^2 + \mu(y) + 1\right)\frac{(b+\exp(z))^3}{a^3\beta^2}(\nabla V(y)g(y))^2$$

$$\dot{V} = \nabla V(y)(f(y)+g(y)k(y)) - C\frac{(b+\exp(z))^3}{a^3\beta^2}\mu^2(y)(\nabla V(y)g(y))^4$$
$$+\nabla V(y)g(y)(\varphi'(y,w)\theta + A'(y,w)d) \quad (4.13)$$
$$-C\left(|A(y,0)|^2 + |\varphi(y,0)|^2 + \mu(y) + 1\right)\frac{(b+\exp(z))^3}{a^3\beta^2}(\nabla V(y)g(y))^2$$

Exploiting (2.3), (2.4), (4.11), (4.12), (4.13), and the facts that $C \geq 1$, $2a\beta < br$, we get for all $(y,w,d,\theta,\delta) \in \mathbb{R}^{n+l+m+p+q}$:

$$\dot{U} \leq -\left(\frac{r}{2} - \frac{a\beta}{b}\right)(R(w)+Q(y))$$
$$+\frac{r}{2}\gamma(|\delta|) + a\frac{|d|^2 + |d|^4 + \left((|\theta|-b-\exp(z))^+\right)^2 + \left((|\theta|-b-\exp(z))^+\right)^4 + 2\beta\Lambda}{b+\exp(z)} \quad (4.14)$$

$$\dot{V} \leq -\frac{r}{2}Q(y) + a\frac{|d|^2 + |d|^4 + \left((|\theta|-b-\exp(z))^+\right)^2 + \left((|\theta|-b-\exp(z))^+\right)^4 + \beta(R(w)+2\Lambda)}{b+\exp(z)} \quad (4.15)$$

Since $V, Q, R, \Phi$ are positive definite and radially unbounded functions and since $2a\beta < br$, there exists $\rho \in K_\infty$, $\tilde{\rho}, \bar{\rho}, M, \tilde{\gamma} \in K_\infty$ such that

$$\left(\frac{r}{2} - \frac{a\beta}{b}\right)(R(w)+Q(y)) \geq \tilde{\rho}(U(y,w)) \quad (4.16)$$

$$\frac{r}{2}Q(y) \geq \rho(V(y)), \quad Q(y) \leq \tilde{\gamma}(V(y)) \quad (4.17)$$

$$M(\Phi(w)) \geq R(w) \geq \bar{\rho}(\Phi(w)) \quad (4.18)$$

Combining the above inequalities with (2.4), (4.14), we obtain for all $(y,w,d,\theta,\delta) \in \mathbb{R}^{n+l+m+p+q}$:

$$\dot{U} \leq -\tilde{\rho}(U(y,w)) + \frac{r}{2}\gamma(|\delta|) + a\frac{|d|^2 + |d|^4 + \left((|\theta|-b-\exp(z))^+\right)^2 + \left((|\theta|-b-\exp(z))^+\right)^4 + 2\beta\Lambda}{b+\exp(z)} \quad (4.19)$$

$$\dot{\Phi} \leq -\bar{\rho}(\Phi(w)) + \bar{\gamma}(|\delta| + V(y)) \quad (4.20)$$

where $\bar{\gamma}(s) := \gamma(s) + \tilde{\gamma}(s)$ for $s \geq 0$.

We next notice that there exist functions $\omega \in KL$, $\zeta \in K_\infty$ with the following property: for every $T \in (0, +\infty]$, $\tilde{\mu} \geq 0$ and for every absolutely continuous function $Y:[0,T) \to \mathbb{R}_+$ for which the following inequality holds for $t \in [0,T)$ a.e.:

$$\dot{Y}(t) \leq -\tilde{\rho}(Y(t)) + \tilde{\mu} \quad (4.21)$$

the following estimate holds for all $t \in [0,T)$:

$$Y(t) \leq \omega(Y(0), t) + \zeta(\tilde{\mu}) \quad (4.22)$$

The existence of functions $\omega \in KL$, $\zeta \in K_\infty$ follows from Lemma 2.14 on page 82 in [6] and by noticing that (4.21) shows that the implication

$$Y(t) \geq \tilde{\rho}^{-1}(2\tilde{\mu}) \Rightarrow \dot{Y}(t) \leq -\tilde{\rho}(Y(t))/2$$

holds for $t \in [0,T)$ a.e..

Similarly, there exist functions $\bar{\omega} \in KL$, $\bar{\zeta} \in K_\infty$ with the following property: for every $T \in (0, +\infty]$, $\tilde{\mu} \geq 0$ and for every absolutely continuous function $Y:[0,T) \to \mathbb{R}_+$ for which the following inequality holds for $t \in [0,T)$ a.e.:

$$\dot{Y}(t) \leq -\bar{\rho}(Y(t)) + \tilde{\mu} \quad (4.23)$$

the following estimate holds for all $t \in [0,T)$:

$$Y(t) \leq \bar{\omega}(Y(0), t) + \bar{\zeta}(\tilde{\mu}) \quad (4.24)$$

Combining (2.14), (4.15), (4.18) we obtain for all $(y, w, d, \theta, \delta) \in \mathbb{R}^{n+l+m+p+q}$:

$$\dot{V} \leq -\rho(V(y))$$
$$+ a\frac{|d|^2 + |d|^4 + \left((|\theta|-b-\exp(z))^+\right)^2 + \left((|\theta|-b-\exp(z))^+\right)^4 + 2\beta\Lambda + \beta M\left(2r^{-1}U(y,w)\right)}{b+\exp(z)} \quad (4.25)$$

By virtue of Lemma 1 there exists a positive non-increasing function $c_\varepsilon : \mathbb{R}_+ \to (0, +\infty)$ such that for every absolutely continuous function $y:[0,T) \to \mathbb{R}_+$ and for every non-increasing function $\alpha \in C^0([0,T]; \mathbb{R}_+)$ for which the differential inequality (2.15) holds for $t \in [0,T)$ a.e., the estimate holds (2.16) for all $t_0 \in [0,T)$ and $t \in [t_0, T)$.

Moreover, define for $s \geq 0$:

$$\tilde{A}(s) := 1 + \max\left(\exp(s), \frac{2aS}{\varepsilon c_\varepsilon(\overline{S})} - b\right), \tag{4.26}$$

$$G(s) := \ln\left(\tilde{A}(s) + \frac{\Gamma \overline{S}}{c_\varepsilon(\overline{S})}\right). \tag{4.27}$$

where $S := 2s^2 + 2s^4 + 2\beta\Lambda + \beta M\left(2r^{-1}\omega(s,0) + 2r^{-1}\zeta(\chi(s,s,s,0))\right)$ and $\overline{S} = s + \rho^{-1}(ab^{-1}S)$. Notice that both $A$ and $G$ are non-decreasing functions.

Let arbitrary $(y_0, w_0, z_0) \in \mathbb{R}^{n+l+1}$, $d \in L^\infty(\mathbb{R}_+; \mathbb{R}^m)$, $\theta \in L^\infty(\mathbb{R}_+; \mathbb{R}^p)$, $\delta \in L^\infty(\mathbb{R}_+; \mathbb{R}^q)$ be given. Using standard theory we conclude that there exists a unique solution $(y(t), w(t), z(t))$ of (2.1), (2.2), (2.7), (2.8) with $(y(0), w(0), z(0)) = (y_0, w_0, z_0)$ defined on a maximal interval $[0, t_{max})$, where $t_{max} \in (0, +\infty]$. For every $t \in [0, t_{max})$, we have from (2.8) that $\dot{z}(t) \geq 0$. Therefore, $z(t)$ is non-decreasing on $[0, t_{max})$ and $z(t) \geq z_0$ for all $t \in [0, t_{max})$. Moreover, the mapping

$$z \mapsto a \frac{|d|^2 + |d|^4 + \left((|\theta| - b - \exp(z))^+\right)^2 + \left((|\theta| - b - \exp(z))^+\right)^4 + 2\beta\Lambda}{b + \exp(z)}$$

is non-increasing. Using (4.19), we conclude that the following inequality holds almost everywhere in $[0, t_{max})$:

$$\frac{d}{dt}(U(y(t), w(t))) \leq -\rho(U(y(t), w(t))) + \chi(\|d\|_\infty, \|\delta\|_\infty, \|\theta\|_\infty, \exp(z_0)) \tag{4.28}$$

Therefore, (4.21), (4.22) imply that estimate (2.11) holds for all $t \in [0, t_{max})$.

It follows from (4.25) and estimate (2.11) that the following inequality holds for $t \in [0, t_{max})$ a.e.:

$$\frac{d}{dt}(V(y(t))) \leq -\rho(V(y(t))) + \frac{aS}{b + \exp(z(t))} \tag{4.29}$$

where

$$\begin{aligned} s &:= U(y_0, w_0) + |z_0| + \|\theta\|_\infty + \|d\|_\infty + \|\delta\|_\infty \\ S &:= 2s^2 + 2s^4 + 2\beta\Lambda + \beta M\left(2r^{-1}\omega(s,0) + 2r^{-1}\zeta(\chi(s,s,s,0))\right) \end{aligned} \tag{4.30}$$

It follows from the fact that $z(t)$ is non-decreasing on $[0, t_{max})$, (4.29), (4.30), (2.14) and (2.16) that the following estimate holds for all $t_0 \in [0, t_{max})$ and $t \in [t_0, t_{max})$:

$$V(y(t)) \leq \min\left(\bar{S}, \bar{S}\exp\left(-c_\varepsilon(\bar{S})(t-t_0)\right) + \frac{\varepsilon}{2} + \frac{aS}{c_\varepsilon(\bar{S})(b+\exp(z(t_0)))}\right) \quad (4.31)$$

where $\bar{S} = s + \rho^{-1}(ab^{-1}S)$.

Next, we show that $z(t)$ is bounded on $[0, t_{max})$. Since $z(t)$ is non-decreasing on $[0, t_{max})$, we have $z(t) \geq z_0$ for all $t \in [0, t_{max})$. We distinguish the following cases:

<u>Case 1:</u> $\exp(z(t)) \leq \tilde{A}(s)$ for all $t \in [0, t_{max})$

In this case, it follows from (4.27) that $z(t) \leq G(s)$, for all $t \in [0, t_{max})$.

<u>Case 2:</u> There exists $T \in (0, t_{max})$ such that $\exp(z(T)) > \tilde{A}(s)$. Since $\exp(z(0)) \leq \exp(s) < \tilde{A}(s)$, (recall (4.26), (4.30)), continuity of $z(t)$ guarantees that there exist $t_0 \in (0, T)$ such that:

$$\exp(z(t_0)) = \tilde{A}(s) \quad (4.32)$$

It follows from (4.26) that $\exp(z(t_0)) \geq \frac{2aS}{\varepsilon c_\varepsilon(\bar{S})} - b$ and consequently, $\frac{aS}{c_\varepsilon(\bar{S})(b+\exp(z(t_0)))} \leq \frac{\varepsilon}{2}$.

Therefore, we obtain from (4.31) for all $t \in [t_0, t_{max})$:

$$V(y(t)) \leq \min\left(\bar{S}, \bar{S}\exp\left(-c_\varepsilon(\bar{S})(t-t_0)\right) + \varepsilon\right) \quad (4.33)$$

Using (2.8) and inequality (4.33), we have for all $t \in [t_0, t_{max})$:

$$\frac{d}{dt}\left(\exp(z(t))\right) \leq \Gamma\bar{S}\exp\left(-c_\varepsilon(\bar{S})(t-t_0)\right) \quad (4.34)$$

Integrating (4.34) and using (4.32), we obtain for all $t \in [t_0, t_{max})$:

$$\exp(z(t)) \leq \exp(z(t_0)) + \frac{\Gamma\bar{S}}{c_\varepsilon(\bar{S})} \leq \tilde{A}(s) + \frac{\Gamma\bar{S}}{c_\varepsilon(\bar{S})} \quad (4.35)$$

Since $z(t)$ is non-decreasing, inequality (4.35) holds for all $t \in [0, t_{max})$.

It follows from (4.35) and definitions (4.26), (4.27), (4.30) that in both cases the following estimates hold for all $t \in [0, t_{max})$:

$$z_0 \leq z(t) \leq G(s) \quad (4.36)$$

$$|z(t)| \leq G(s) \quad (4.37)$$

Estimate (4.37) implies that $z(t)$ is bounded on $[0, t_{max})$.

Definition (2.14), estimate (2.11) and the fact that $V, \Phi$ are positive definite and radially unbounded functions guarantee that the component $(y(t), w(t))$ of the solution $(y(t), w(t), z(t))$ of (2.1), (2.2), (2.7), (2.8) with $(y(0), w(0), z(0)) = (y_0, w_0, z_0)$ is bounded on $[0, t_{max})$.

Consequently, $t_{max} = +\infty$. Moreover, estimate (2.11) holds for all $t \geq 0$. Estimate (2.11), equation (2.7) and the fact that $V, \Phi$ are positive definite and radially unbounded functions allow us to conclude that $y \in L^\infty(\mathbb{R}_+; \mathbb{R}^n)$, $w \in L^\infty(\mathbb{R}_+; \mathbb{R}^l)$, $z \in L^\infty(\mathbb{R}_+; \mathbb{R})$ and $u \in L^\infty(\mathbb{R}_+; \mathbb{R})$.

Since $z(t)$ is non-decreasing and bounded from above, $\lim_{\tau \to +\infty}(z(\tau))$ exists. Thus, we obtain (2.12) from (4.36) with $B(y_0, w_0, z_0, \|d\|_\infty, \|\delta\|_\infty, \|\theta\|_\infty) = \tilde{G}(U(y_0, w_0) + |z_0| + \|\theta\|_\infty + \|d\|_\infty + \|\delta\|_\infty)$ and $\tilde{G}(s) := \int_0^{s+1} G(l) dl$ for $s \geq 0$ (notice that $\tilde{G}: \mathbb{R}_+ \to \mathbb{R}_+$ is a non-decreasing, continuous function that satisfies $\tilde{G}(s) \geq G(s)$ for all $s \geq 0$).

We next show estimate (2.9). Estimate (4.36) and the fact that $z(t)$ is non-decreasing guarantees that the function $z(t)$ has a finite limit as $t \to +\infty$. This implies that the function $\exp(z(t))$ has a finite limit as $t \to +\infty$. Moreover, the facts that $d \in L^\infty(\mathbb{R}_+; \mathbb{R}^m)$, $\theta \in L^\infty(\mathbb{R}_+; \mathbb{R}^p)$, $y \in L^\infty(\mathbb{R}_+; \mathbb{R}^n)$, $u \in L^\infty(\mathbb{R}_+)$, $w \in L^\infty(\mathbb{R}_+; \mathbb{R}^l)$ and (2.1) imply that $\frac{d}{dt}(V(y(t)))$ is of class $L^\infty(\mathbb{R}_+)$. It follows that the function $\frac{d}{dt}(\exp(z(t))) = \Gamma(V(y(t)) - \varepsilon)^+$ is uniformly continuous, in addition to $\int_0^{+\infty} \frac{d}{dt}(\exp(z(t))) dt = \lim_{t \to +\infty}(\exp(z(t))) - \exp(z_0) \leq \exp(G(s)) - \exp(z_0) < +\infty$. From Barbălat's Lemma (see [12, 7]), we have:

$$\lim_{t \to +\infty}\left(\frac{d}{dt}(\exp(z(t)))\right) = \lim_{t \to +\infty}\left(\Gamma(V(y(t)) - \varepsilon)^+\right) = 0 \qquad (4.38)$$

Therefore, estimate (2.9) holds.

Exploiting (4.20), (4.23), (4.24) we get for all $T > 0$ and $t \geq T$:

$$\Phi(w(t)) \leq \bar{\omega}(\Phi(w(T)), t - T) + \bar{\zeta}\left(\bar{\gamma}\left(\sup_{s \geq T}(|\delta(s)| + V(y(s)))\right)\right) \qquad (4.39)$$

Using estimate (2.11) and definition (2.14) in conjunction with (4.39), we get for all $T > 0$ and $t \geq T$:

$$\Phi(w(t)) \leq \bar{\zeta}\left(\bar{\gamma}\left(\sup_{s \geq T}(|\delta(s)| + V(y(s)))\right)\right)$$
$$+ \bar{\omega}\left(2r^{-1}\omega(U(y_0, w_0), 0) + 2r^{-1}\zeta\left(\chi(\|d\|_\infty, \|\delta\|_\infty, \|\theta\|_\infty, \exp(z_0))\right), t - T\right) \quad (4.40)$$

Consequently, we obtain from (4.40) for all $T > 0$:

$$\limsup_{t \to +\infty}(\Phi(w(t))) \leq \bar{\zeta}\left(\bar{\gamma}\left(\sup_{s \geq T}(|\delta(s)| + V(y(s)))\right)\right) \quad (4.41)$$

Next, we show estimate (2.10). Let arbitrary $\epsilon > 0$ be given. Then there exists $T > 0$ such that

$$|\delta(t)| \leq \limsup_{\tau \to +\infty}(|\delta(\tau)|) + \epsilon \text{ for } t \geq T \text{ a.e.}$$

$$V(y(t)) \leq \limsup_{\tau \to +\infty}(V(y(\tau))) + \epsilon \text{ for all } t \geq T$$

Therefore, we get from (4.41):

$$\limsup_{t \to +\infty}(\Phi(w(t))) \leq \bar{\zeta}\left(\bar{\gamma}\left(\limsup_{\tau \to +\infty}(|\delta(\tau)|) + \limsup_{\tau \to +\infty}(V(y(\tau))) + 2\epsilon\right)\right) \quad (4.42)$$

Estimate (2.10) with $\tilde{\zeta} = \bar{\zeta} \circ \bar{\gamma}$ is a direct consequence of (4.42), (2.9) and the fact that $\epsilon > 0$ is arbitrary. Since $\bar{\gamma}(s) := \gamma(s) + \tilde{\gamma}(s)$ and $\tilde{\gamma}, \bar{\zeta} \in K_\infty$ it follows that $\tilde{\zeta}(0) = 0$ when $\gamma(0) = 0$.

Finally, combining (4.19), (4.22) and the semigroup property we obtain:

$$\limsup_{t \to +\infty}\left(\zeta^{-1}(U(y(t), w(t)))\right) \leq \frac{r}{2}\gamma\left(\limsup_{t \to +\infty}(|\delta(t)|)\right)$$
$$+ a\frac{\left(\limsup_{t \to +\infty}(|d(t)|)\right)^2 + \left(\limsup_{t \to +\infty}(|d(t)|)\right)^4 + 2\beta\Lambda}{b + \exp\left(\lim_{t \to +\infty}(z(t))\right)}$$
$$+ a\frac{\left(\left(\limsup_{t \to +\infty}(|\theta(t)|) - b - \exp\left(\lim_{t \to +\infty}(z(t))\right)\right)^+\right)^2 + \left(\left(\limsup_{t \to +\infty}(|\theta(t)|) - b - \exp\left(\lim_{t \to +\infty}(z(t))\right)\right)^+\right)^4}{b + \exp\left(\lim_{t \to +\infty}(z(t))\right)} \quad (4.43)$$

The asymptotic estimate (4.43) shows that if $\gamma(0) = \Lambda = 0$, $\lim_{t \to +\infty}(d(t)) = 0$, $\lim_{t \to +\infty}(\delta(t)) = 0$ and $\theta$ is constant then we have:

$$\limsup_{t \to +\infty}\left(\zeta^{-1}(U(y(t), w(t)))\right)$$
$$\leq a\frac{\left(\left(|\theta| - b - \exp\left(\lim_{t \to +\infty}(z(t))\right)\right)^+\right)^2 + \left(\left(|\theta| - b - \exp\left(\lim_{t \to +\infty}(z(t))\right)\right)^+\right)^4}{b + \exp\left(\lim_{t \to +\infty}(z(t))\right)} \quad (4.44)$$

If $\limsup\limits_{t\to+\infty}\left(\zeta^{-1}(U(y(t),w(t)))\right)>0$ then $\lim\limits_{s\to+\infty}(z(s))<\ln(|\theta|-b)$. Otherwise, $\limsup\limits_{t\to+\infty}\left(\zeta^{-1}(U(y(t),w(t)))\right)=0$ which implies that $\lim\limits_{t\to+\infty}(|y(t)|)=\lim\limits_{t\to+\infty}(|w(t)|)=0$.

The proof is complete. ◁

We finally provide the proof of Theorem 2. The proof of Theorem 2 is based on three facts:

(i) For every $a\in K_\infty$ there exists $\bar{A}\in C^\infty(\mathbb{R}^n)$ being positive definite and radially unbounded and $\bar{\Lambda}>0$ such that $a(|x|)\le\bar{\Lambda}+\bar{A}(x)$ for all $x\in\mathbb{R}^n$,

(ii) For every $a\in K_\infty$ there exists $\hat{A}\in C^\infty(\mathbb{R}_+;\mathbb{R}_+)$ being non-decreasing such that $a(s)\le\hat{A}(s)$ for all $s\ge 0$,

(iii) For every $\tilde{A}\in C^0(\mathbb{R}^n)$ with $\tilde{A}(0)=0$ there exists $\tilde{a}\in K_\infty$ such that $|\tilde{A}(x)|\le\tilde{a}(|x|)$ for all $x\in\mathbb{R}^n$.

**Proof of Theorem 2:** We show that assumption (A) is valid for appropriate smooth mappings $V,Q:\mathbb{R}^n\to\mathbb{R}_+$, $\Phi,R:\mathbb{R}^l\to\mathbb{R}_+$, $\mu:\mathbb{R}^n\to(0,+\infty)$, $V,Q,R,\Phi$ being positive definite and radially unbounded, appropriate constants $r,\Lambda>0$ and an appropriate non-decreasing function $\gamma\in C^0(\mathbb{R}_+;\mathbb{R}_+)$. We construct all the above functions and constants step-by-step.

<u>Step 1:</u> Construction of $R:\mathbb{R}^l\to\mathbb{R}_+$ and $\mu:\mathbb{R}^n\to(0,+\infty)$.

Clearly, there exists $a\in K_\infty$ such that

$$|\varphi(y,w)-\varphi(y,0)|^2+|A(y,w)-A(y,0)|^2\le a(|(y,w)|) \tag{4.45}$$

for all $(y,w)\in\mathbb{R}^{n+l}$. Therefore, we get from (4.45) for all $(y,w)\in\mathbb{R}^{n+l}$:

$$\begin{aligned}|\varphi(y,w)-\varphi(y,0)|^2+|A(y,w)-A(y,0)|^2 &\le a(|y|+|w|)\\ &\le a(2|y|)+a(2|w|)\le(1+a(2|y|))(1+a(2|w|))\end{aligned} \tag{4.46}$$

Clearly, there exist $\bar{A}\in C^\infty(\mathbb{R}^n)$, $R\in C^\infty(\mathbb{R}^l)$ being positive definite and radially unbounded and $\Lambda_1>0$ such that $a(2|w|)\le\Lambda_1+R(w)$ for all $w\in\mathbb{R}^l$ and $a(2|y|)\le\Lambda_1+\bar{A}(y)$ for all $y\in\mathbb{R}^n$. Therefore, we obtain from (4.46) for all $(y,w)\in\mathbb{R}^{n+l}$:

$$|\varphi(y,w)-\varphi(y,0)|^2+|A(y,w)-A(y,0)|^2\le(1+\Lambda_1+\bar{A}(y))(1+\Lambda_1+R(w)) \tag{4.47}$$

Moreover, there exists $\bar{a}\in K_\infty$ such that the following inequality holds for all $y\in\mathbb{R}^n$:

$$\frac{|\varphi(y,0)|^2}{1+\Lambda_1+\bar{A}(y)} \leq \bar{a}(|y|) \tag{4.48}$$

Clearly, there exist $\tilde{A} \in C^\infty(\mathbb{R}^n)$ being positive definite and radially unbounded and $\Lambda_2 > 0$ such that $\bar{a}(|y|) \leq \Lambda_2 + \tilde{A}(y)$ for all $y \in \mathbb{R}^n$. Therefore, we obtain from (4.48) for all $y \in \mathbb{R}^n$:

$$\frac{|\varphi(y,0)|^2}{1+\Lambda_1+\bar{A}(y)} \leq \Lambda_2 + \tilde{A}(y) \tag{4.49}$$

By defining for all $y \in \mathbb{R}^n$

$$\mu(y) := 1 + \Lambda_1 + \bar{A}(y) \tag{4.50}$$

we obtain from (4.47) and (4.49) that the following inequalities hold for all $(y,w) \in \mathbb{R}^{n+l}$:

$$\begin{aligned}|\varphi(y,0)|^2 &\leq \mu(y)\big(\Lambda_2 + \tilde{A}(y)\big) \\ |\varphi(y,w)-\varphi(y,0)|^2 + |A(y,w)-A(y,0)|^2 &\leq \mu(y)\big(1+\Lambda_1+R(w)\big)\end{aligned} \tag{4.51}$$

<u>Step 2:</u> Construction of $\Phi : \mathbb{R}^l \to \mathbb{R}_+$.

By virtue of Theorem 2.63 in [19] and since subsystem (2.2) with inputs $(y,\delta) \in \mathbb{R}^n \times \mathbb{R}^q$ satisfies the ISS property, there exist a smooth mapping $\tilde{\Phi} : \mathbb{R}^l \to \mathbb{R}_+$ being positive definite and radially unbounded and $\tilde{\gamma}, \tilde{a} \in K_\infty$ such that the following implication holds for all $(y,w,\delta) \in \mathbb{R}^n \times \mathbb{R}^l \times \mathbb{R}^q$:

$$\tilde{\Phi}(w) \geq \tilde{\gamma}(|(y,\delta)|) \Rightarrow \nabla \tilde{\Phi}(w) h(y,w,\delta) \leq -\tilde{a}\big(\tilde{\Phi}(w)\big) \tag{4.52}$$

Let $\Lambda_3 > 0$ be a given constant. We define the function $\tilde{g} : \mathbb{R}^l \to \mathbb{R}_+$ by means of the following formulas:

$$\tilde{g}(w) = \left(\frac{R(w)-\Lambda_3}{\tilde{a}\big(\tilde{\Phi}(w)\big)}\right)^+ \text{ for } w \neq 0 \text{ and } \tilde{g}(0) = 0 \tag{4.53}$$

Definition (4.53) shows that $\tilde{g}$ is a non-negative, continuous function. By virtue of Lemma 2.4 in [6] there exists $\bar{g} \in K_\infty$ such that the following inequality holds for all $s \geq 0$:

$$\max\{\tilde{g}(w) : \tilde{\Phi}(w) \leq s, w \in \mathbb{R}^l\} \leq \bar{g}(s) \tag{4.54}$$

Moreover, there exists $\hat{A} \in C^\infty(\mathbb{R}_+;\mathbb{R}_+)$ being non-decreasing such that $\bar{g}(s) \leq \hat{A}(s)$ for all $s \geq 0$. Combining (4.53), (4.54) we get for all $w \in \mathbb{R}^l$:

$$\hat{A}\big(\tilde{\Phi}(w)\big)\tilde{a}\big(\tilde{\Phi}(w)\big) \geq R(w) - \Lambda_3 \tag{4.55}$$

Define the functions

$$\tilde{b}(s) := s + \int_0^s \hat{A}(\eta)d\eta, \quad \Phi(w) := \tilde{b}\left(\tilde{\Phi}(w)\right) \tag{4.56}$$

for all $s \geq 0$, $w \in \mathbb{R}^l$. Notice that $\tilde{b} \in K_\infty$ and consequently, $\Phi : \mathbb{R}^l \to \mathbb{R}_+$ defined by (4.56) is a smooth, radially unbounded, positive definite function. Furthermore, implication (4.52) and (4.55) imply that the following implication holds for all $(y, w, \delta) \in \mathbb{R}^n \times \mathbb{R}^l \times \mathbb{R}^q$:

$$\Phi(w) \geq \tilde{b}\left(\tilde{\gamma}(|(y,\delta)|)\right) \Rightarrow \nabla \Phi(w) h(y, w, \delta) \leq -R(w) + \Lambda_3 \tag{4.57}$$

Next define for all $s \geq 0$:

$$\rho(s) := \max\left\{ \nabla \Phi(w) h(y, w, \delta) + R(w) : \Phi(w) \leq \tilde{b}\left(\tilde{\gamma}(|(y,\delta)|)\right), |(y,\delta)| \leq s \right\} \tag{4.58}$$

The fact that $\nabla \Phi(0) = 0$ (a consequence of Fermat's theorem and the fact that $\Phi : \mathbb{R}^l \to \mathbb{R}_+$ is a smooth, positive definite function) implies that $\rho(s)$, as defined by (4.58), is continuous at $s = 0$ with $\rho(0) = 0$. Consequently, since $\rho(s)$, as defined by (4.58), is non-decreasing, there exists $\bar{\gamma} \in K_\infty$ such that $\rho(s) \leq \bar{\gamma}(s)$ for all $s \geq 0$. Combining (4.57) and (4.58) we get for all $(y, w, \delta) \in \mathbb{R}^n \times \mathbb{R}^l \times \mathbb{R}^q$:

$$\nabla \Phi(w) h(y, w, \delta) \leq -R(w) + \Lambda_3 + \bar{\gamma}(|(y,\delta)|)$$
$$\leq -R(w) + \Lambda_3 + \bar{\gamma}(2|\delta|) + \bar{\gamma}(2|y|)$$

Clearly, there exist $\bar{B} \in C^\infty(\mathbb{R}^n)$ being positive definite and radially unbounded and $\Lambda_4 > 0$ such that $\bar{\gamma}(2|y|) \leq \Lambda_4 + \bar{B}(y)$ for all $y \in \mathbb{R}^n$. Thus, we obtain for all $(y, w, \delta) \in \mathbb{R}^n \times \mathbb{R}^l \times \mathbb{R}^q$:

$$\nabla \Phi(w) h(y, w, \delta) \leq -R(w) + \Lambda_3 + \Lambda_4 + \bar{\gamma}(2|\delta|) + \bar{B}(y) \tag{4.59}$$

<u>Step 3:</u> Construction of $Q : \mathbb{R}^n \to \mathbb{R}_+$, $V : \mathbb{R}^n \to \mathbb{R}_+$ and $r > 0$.

Using the results in [16] and since $0 \in \mathbb{R}^n$ is globally asymptotically stable for the closed-loop system (2.1) with $\theta = 0$, $d = 0$ and $u = k(y)$, there exist a smooth mapping $\tilde{V} : \mathbb{R}^n \to \mathbb{R}_+$ being positive definite and radially unbounded and $\bar{b} \in K_\infty$ such that the following inequality holds for all $y \in \mathbb{R}^n$:

$$\nabla \tilde{V}(y)\left(f(y) + g(y)k(y)\right) \leq -\bar{b}\left(\tilde{V}(y)\right) \tag{4.60}$$

We also define for all $y \in \mathbb{R}^n$

$$\tilde{Q}(y) := \tilde{A}(y) + \bar{B}(y) \tag{4.61}$$

Notice that $\tilde{Q} : \mathbb{R}^n \to \mathbb{R}_+$ is a smooth, radially unbounded, positive definite function.

We define the function $\tilde{G} : \mathbb{R}^n \to \mathbb{R}_+$ by means of the following formulas:

$$\tilde{G}(y) = \left( \frac{\tilde{Q}(y) - \Lambda_3}{\bar{b}(\tilde{V}(y))} \right)^+ \text{ for } y \neq 0 \text{ and } \tilde{G}(0) = 0 \tag{4.62}$$

Definition (4.62) shows that $\tilde{G}$ is a non-negative, continuous function. By virtue of Lemma 2.4 in [6] there exists $\bar{G} \in K_\infty$ such that the following inequality holds for all $s \geq 0$:

$$\max\left\{ \tilde{G}(y) : \tilde{V}(y) \leq s, y \in \mathbb{R}^n \right\} \leq \bar{G}(s) \tag{4.63}$$

Moreover, there exists $\hat{B} \in C^\infty(\mathbb{R}_+; \mathbb{R}_+)$ being non-decreasing such that $\bar{G}(s) \leq \hat{B}(s)$ for all $s \geq 0$. Combining (4.62), (4.63) we get for all $y \in \mathbb{R}^n$:

$$\hat{B}(\tilde{V}(y))\bar{b}(\tilde{V}(y)) \geq \tilde{Q}(y) - \Lambda_3 \tag{4.64}$$

Define the functions

$$\tilde{B}(s) := s + \int_0^s \hat{B}(\eta) d\eta, \; V(y) := \tilde{B}(\tilde{V}(y)), \; Q(y) := \left(1 + \hat{B}(\tilde{V}(y))\right)\bar{b}(\tilde{V}(y)) \tag{4.65}$$

Definition (4.65) shows that $V, Q : \mathbb{R}^n \to \mathbb{R}_+$ are both smooth, positive definite and radially unbounded functions for which (2.3) holds for all $y \in \mathbb{R}^n$ with $r = 1$.

<u>Step 4:</u> Construction of $\gamma \in C^0(\mathbb{R}_+; \mathbb{R}_+)$ and $\Lambda > 0$.

Using (4.64), (4.65), (4.61) and (4.51) we conclude that (2.5) and (2.6) hold with $\Lambda := 1 + \Lambda_1 + \Lambda_2 + \Lambda_3$.

Using (4.64), (4.65), (4.61) and (4.59) we conclude that (2.4) holds with $\gamma(s) := 2\Lambda_3 + \Lambda_4 + \bar{\gamma}(2s)$. The proof is complete. ◁

## 5. Concluding Remarks

The important (and rarely achieved) combination of robustness properties that are guaranteed by the DADS controller for the closed-loop system in conjunction with its simplicity justify our focus to a special class of systems (systems with matched unmodeled dynamics). However, a detailed study is needed for the characterization of the class of nonlinear systems for which these robustness properties can be guaranteed by adaptive control schemes.

It is clear that the procedure described in Section 3 can be applied -in exactly the same way- to different PDEs with different boundary conditions. Moreover, nonlinear functions $K : [0,1] \times \mathbb{R} \to \mathbb{R}$ that satisfy other bounds than (3.2) can also be studied. Therefore, the partial-state DADS feedback design is not limited to finite-dimensional systems and many other systems with infinite-dimensional uncertain dynamics can be studied.